\documentstyle{amsppt}
\magnification=\magstephalf
\NoRunningHeads
\NoBlackBoxes

\title
A factorization constant for $l^n_p$, $0 < p < 1$
\endtitle

\author
N.T. Peck
\endauthor

\endtopmatter
\def\EOP{\unskip\nobreak\hfill$\blacksquare$}

\document
\baselineskip=24pt

In this paper we are concerned with factoring the identity 
operator on an $n$-dimensional quasi-normed space $X$ through a 
space $\ell^K_{\infty}$.

We seek a good lower bound 
${\underline\lambda}(x)$
for $\parallel P \parallel\,\, \parallel T \parallel$ 
over all factorizations $Id_X = PT,$
with $T:X \longrightarrow \ell_{\infty}^K,
\ \ P:\ell_{\infty}^K \longrightarrow X .$
When $X$ is $\ell^n_p$, $1 \leq p$, the 
constant is known: see \cite{5, Theorem 32.9} and the references 
given for that theorem.

For $p < 1$, we will obtain the lower estimate $\underline 
\lambda (l^n_p) \geq Cn^{\frac1p - \frac12}(\log n)^{-\frac12}$. 
(A $T$ and $P$ with $\parallel P \parallel\, \, \parallel T \parallel 
\leq Cn^{\frac1p - \frac12}$ are easily obtained.)

Throughout, $C$ denotes a constant, which may vary from one 
occurrence to the next, but which is independent of $n$.

We thank Y. Gordon for valuable conversations.

\proclaim{Lemma 1}
Let $(r_i)$, $1 \leq i \leq n$, be the first $n$ Rademacher 
functions and let $\alpha_1, \dots\alpha_n$ be real.  Letting $m$ 
denote Lebesque measure on $(0,1)$, we have
$$
\align
m \Biggl\{ \left| \sum_{i=1}^n \alpha_i r_i\right | &> \alpha 
\sqrt{\log n \sum \alpha_i^2} \Biggr\} \\
&\leq n^{-C\alpha^2},                                        
\endalign
$$
for any positive $\alpha$.
\endproclaim

\demo{Proof}
This is well known; for completeness, we sketch a proof, 
following a suggestion of R. Kaufman.

We can assume $\sum_{i=1}^n\alpha_i^2 = 1$.  Put  $f = 
\sum_{i=1}^n\alpha_ir_i$.  By Khintchine's inequality, for some 
constant $C$ and all $p \geq 1$, $\|f\|_p \leq Cp^{\frac12}$.  
>From this, $m\{|f| > \lambda\} \leq (C/\lambda)^pp^{\frac p2} = 
K^pp^{\frac p2}$ with $K = C/\lambda$.

Now minimize $K^pp^{\frac p2}$ in $p$.  At the minimizer, we find 
$p = K^{-2}e$, from which $K^pp^{\frac p2} = exp(-
\lambda^2/eC^2)$.  (Note that $p > 1$ for $\lambda > C.)$  
Finally, put $\lambda = \alpha\sqrt{\log n}$ to get the 
conclusion.
\EOP

The space $\ell_\infty$ is not of type 2;  but the conclusion of the 
next lemma will suffice for our purposes.
\enddemo

\proclaim{Lemma 2}
Let $Y$ be an $n$-dimensional subspace of $L_1(0,1)$ and let 
$f_1\dots f_n$ be elements of $Y^\ast$, of norm at most $1$.  Then 
for some $\overline s$ in $(0,1)$,
$$
\sup_{\|y\| \leq 1} \left| \sum^n_{i=1} r_i(\overline s)f_i(y)\right|
\leq C\sqrt{n \log n}.
$$
\endproclaim
          
\demo{Proof}
For any $0 < \epsilon < \frac12$, a result of Schechtman \cite{3} 
implies that there are  an $N \leq C\epsilon^{-2} \log(\epsilon^{-
1})n^2$ and an isomorphism $U : Y \overset 
\text{into}\to\longrightarrow \ell^N_1$ such that $\|U\|\,\,\|U^{-
1}\| \leq 1 + \epsilon$. See also the results of
Bourgain--Lindenstrauss--Milman \cite{1} and Talagrand \cite{4}.  In 
particular, taking $\epsilon = \frac14$, say, we obtain the 
corresponding $U$; we can assume $\|U\| = 1$, so that $\|U^{-1}\| 
\leq \frac54$ and $N \leq Cn^2$ (after changing $C$) $\leq n^3$, 
if $n$ is sufficiently large.

Let $e_1,\dots e_N$ be the unit basis vectors in $\ell^N_1$.  For 
each $i$ let $\Phi_i$ be a Hahn--Banach extension to $\ell^N_1$ of 
$f_iU^{-1}$ on $U(Y)$, with $\|\Phi_i\| \leq \frac54$; then for 
each $j$, $\Phi_i(e_j)\| \leq \frac54$.

Now fix $\alpha$ with $C\alpha^2 > 3$, where $C$ is the 
constant in the conclusion of lemma $1$; then $n^3 n^{-C\alpha^2} 
\leq \frac14$ if $n$ is sufficiently large.  Since 
$$
\align
&m \Biggl\{ \left| \sum \phi_i(e_j)r_i\right | > \frac54 \alpha 
\sqrt{n\log n} \Biggr\} \\               
&\leq m \Biggl\{ \left| \sum \phi_i(e_j)r_i\right | > \alpha 
\sqrt{\log n \sum \phi^2_i(e_j)} \Biggr\} \\
&\leq n^{-C\alpha^2},\qquad \text{for each}\quad j,
\endalign
$$
there is a set $A$, $m(A) > \frac34$, such that
$\left| \sum^n_{i=1}\phi_i(e_j)r_i (\overline s)\right | 
\leq \frac54 \alpha\sqrt{n \log n}$ for each $j$ and each 
$\overline s$ in $A$.

Now if $y \in Y$ and $\| y \| \leq 1$, $\|Uy \| \leq 1$.  Write 
$Uy = \sum_{j=1}^N \alpha_j e_j$, $\sum_{j=1}^N |\alpha_j| \leq 
1$; applying the above inequality to each $j$ and recalling that 
$f_i = f_i U^{-1}U$ on $Y$, we have 
$$
\align
\left| \sum_{i=1}^n f_i(y)r_i(\overline s)\right| &\leq \frac54 
\alpha \sqrt{n \log n}\\
&= C\sqrt{n \log n}
\endalign
$$
for each $\overline s$ in $A$.
\EOP
\enddemo
\example{Notation}
Let $\Cal A$ be an algebra of measurable subsets of $(0,1)$. For 
$0 < p \leq \infty$, $L_p(\Cal A)$ is the space of functions in 
$L_p(0,1)$ which are $\Cal A$-measurable.  For ease of argument, 
we deal with an $L_\infty(\Cal A)$ with ``homogeneous'' $\Cal A$, 
rather than $\ell^K_\infty$.
\endexample
\proclaim{Theorem}
Let $\Cal A$ be a finite subalgebra of measurable subsets of $(0,1)$
containing the dyadic intervals $\left( \frac{j - 1}{2^n}, 
\frac{j}{2^n} \right)$, $1 \leq j \leq 2^n$, 
and assume the atoms of $\Cal A$ all 
have the same measure. Let $X$ be an $n$-dimensional vector 
space.  Let $T : X \longrightarrow L_\infty (\Cal A)$ be a linear 
map; and let $(e_i)^n_{i=1}$ be elements of $X$ such that 
$\|T(e_i)\|_\infty \leq 1$, $1 \leq i \leq n$.  Let $P : 
L_\infty(\Cal A) \longrightarrow X$ be a linear operator such 
that $PT = Id_X$.  Then there is $w$ in $L_\infty(\Cal A)$ with 
$\|w\|_\infty \leq C\sqrt{n \log n}$ such that $P(w) = \sum_{i=1}^n 
r_i(\overline s)e_i$, for some $\overline s$ in $(0,1)$.
\endproclaim

\demo{Proof}
Let ${(I_j)}_{j=1}^k$ be the atoms of $\Cal A$, and let
$z_1,\dots z_{k-n}$ be a basis for $\ker P$. Define a $(k-n) \times k$ 
matrix by $z_{i,j} =$ constant value of $z_i$ on the atom $I_j$.

Now row--reduce the matrix $(z_{i,j})$.  In $k - n$ of the 
columns there will be one $1$ with all other entries $0$; denote 
the atoms corresponding to the $n$ remaining ``distinguished'' 
columns by $I_{s_1},\dots I_{s_n}$.  Enlarge the matrix 
$(z_{i,j})$ to a $k \times k$ matrix $(y_{i,j})$ by adding $n$ 
rows of zeros in each of rows $s_1$ through $s_n$.

We can obviously regard $y_{i,j}$ as an $\Cal A \times
\Cal A$-measurable function $y(s,t)$ on $(0,1) \times (0,1)$, which 
satisfies these properties:
\roster
\item"{(1)}" if $s \notin \underset j\to U I_{s_j}$ and if $s$ and 
$t$ are in the same atom, $y(s,t) = 1$;                               
\item"{(2)}" if $s \notin \underset j\to U I_{s_j}$ and if $s$ and 
$t$ are in different atoms, $y(s,t) = 0$;                       
\item"{(3)}" if $t \notin \underset j\to U I_{s_j}$,
$y(s,t) = 0$ for all $s$;
\item"{(4)}" if we define $y_t$ on $(0,1)$ by $y_t(s) = y(s,t)$, 
then $y_t \in \ker P$; 
\item"{(5)}" a function $f$ in $L_\infty(\Cal A)$ is in $\ker P$ 
if and only if there is a function $\beta(t)$ in $L_\infty(\Cal 
A)$ so that $f(s) = \int y_t(s)\beta(t)dt$ for all $s$.
\endroster

Properties $(1)-(5)$ are evident from the description of $\ker P$ 
and the properties of a row-reduced matrix.

Let $s^\ast_i$ be a point of the atom $I_{s_i}$, $1 \leq i \leq n$, 
and let $g_i = T(e_i)$, $1 \leq i \leq n$; then $\|g_i\|_\infty 
\leq 1$.

Now define 
$$
R(s,t) = \sum_{i=1}^n r_i(s)g_i(t),
$$                  
and for a function $y(t)$, define
$$
\align
\psi_s(y) &= \int R(s,t)y(t)dt\\
&=\sum_{i=1}^n (\int y(t)g_i(t)dt)r_i(s).
\endalign
$$

Let $Y$ be the span of the functions $y_t(s^\ast_i)$, $1 \leq i \leq 
n$, regarded as functions of $t$.  Since $y \longrightarrow \int 
y(t) g_i(t)dt$ is of norm at most 1, Lemma 2 implies that there 
is a set $A$ of measure $> \frac34$ such that for any $\overline 
s$ in $A$, 
$$     
|\psi_{\overline s}(y)| \leq C\sqrt{n \log n}\, \|y\|_1
$$
for all $y$ in $Y$.

Take any norm-preserving extension of $\psi_s$ to all of $L_1(\, 
,dt)$; then there is a function $h_{\overline s}(t)$ in 
$L_\infty(.,dt)$ with $\|h_{\overline s}\|_\infty \leq C\sqrt{n 
\log n}$ and such that
$$
\psi_{\overline s}(y) = \int y(t)h_{\overline s}(t)dt
$$
for all $y$ in $Y$.  Restating this, 
$$
\int (R_{\overline s}(t) - h_{\overline s}(t))y(t) dt = 0
\tag6 
$$
for all $y$ in $Y$.  Now put
$$
\alpha_{\overline s}(s) = \frac1m \int (R_{\overline s}(t) - 
h_{\overline s}(t))y_t(s) dt,
$$
where $m$ is the measure of an atom of $\Cal A$.  Then 
$\alpha_{\overline s} \in \ker P$ by property (5), 
$\alpha_{\overline s}(s) = R_{\overline s}(s) - h_{\overline 
s}(s)$, for $s$ not in $\bigcup_{i=1}^n I_{s_i}$, by (1) and (2), 
and $\alpha_{\overline s}(s_i) = 0$, $1 \leq i \leq n$, by (6).

The set $A$ in the proof of Lemma 2 has measure $\geq \frac34$.  
vvApplying Lemma 1 again, we can require that
$$
|R_{\overline s}(s^\ast_i)| \leq C\sqrt{n \log n} 
$$
for each $i$, $1 \leq i \leq n$, for all $\overline s$ in a set 
$B$ with $m(B) \geq \frac34$.  Thus $A \cap B$ has positive 
measure, so choose $\overline s$ in $A \cap B$.

Now, to finish the proof, set $w(s) = R_{\overline s}(s) - 
\alpha_{\overline s}(s)$.  If $s \notin \bigcup_{i=1}^n I_{s_i}$, 
$w(s) = R_{\overline s}(s) - \alpha_{\overline s}(s) = 
h_{\overline s}(s)$, so $|w(s)| \leq C\sqrt{n \log n}$.  Also, 
$|w(s^\ast_i)| = |R_{\overline s}(s^\ast_i)| \leq C\sqrt{n \log n}$; thus 
$\|w\|_\infty \leq C\sqrt{n \log n}$.  Finally, $Pw = 
\sum_{i=1}^n r_i(\overline s)e_i$ since $\alpha_{\overline s} \in 
\ker P$.  This completes the proof. \EOP
\enddemo

\proclaim{Corollary 1}
Let $X$ be $n$-dimensional, and let $X \overset T\to\longrightarrow
\ell_{\infty}^K \overset P\to
\longrightarrow X$  be a factorization of $Id_X$ through
$\ell_{\infty}^K$. Then $\Vert P\Vert\ \Vert T\Vert \geq 
Cs_n(n \log n)^{\frac12}$, where
$$
s_n = \sup\Sb \|Te_i\| \leq 1\\1 \leq i \leq n \endSb \inf_{\pm 1}
\|\sum_{i=1}^n \pm e_i\|.                                  
$$

In particular, if $\|T\| \leq 1$, $\|P\| \geq Cb_n(n \log n)^{-
\frac12}$, where $b_n =  \displaystyle{\sup\limits\Sb\|e_i\| \leq 
1\\1 \leq i 
\leq n\endSb}\inf\limits_{\pm 1}
\|\sum_{i=1}^n \pm e_i\|$.  (See \cite{2}.)
\endproclaim

\demo{Proof} For $K\geq n\ ,$ this is an immediate consequence of the Theorem.
Assume now that $K < n\ .$ Let $j=n-K\ ,$ and define $\tilde T\ :X
\longrightarrow\ell_{\infty}^n=\ell_{\infty}^K\otimes
\ell_{\infty}^j$ by $\tilde T(x)= (T(x),0)\ .$ Define
$\tilde P : \ell_{\infty}^n\longrightarrow X$ by $\tilde P(w,z)=
P(w) .$ Note that $\tilde P\tilde T= Id_X\ $ and that $\Vert\tilde P\Vert=
\Vert P\Vert\ , \Vert\tilde T\Vert=\Vert T\Vert\ .$ The result now follows
from the Theorem. \EOP
\enddemo

\proclaim{Corollary 2}
Suppose $\ell_p^n\overset T\to\longrightarrow\ell_{\infty}^K
\overset P\to\longrightarrow \ell^n_p$ is any factorization of 
the identity on $\ell^n_p$ through $\ell^K_\infty$, $0 < p \leq 1$, 
with $\|T\| = 1$.  Then $\|P\| \geq Cn^{\frac1p - \frac12}(\log 
n)^{-\frac12}$. 
\endproclaim

\demo{Proof}
Take $X = \ell^n_p$, $e_i$ the usual $i$th basis vector in 
$\ell^n_p$, $1 \leq i \leq n$, for $0 < p \leq 1$. Now apply
Corollary 1. \EOP
\enddemo

\example
{Remark}
For $0 < p < 1$, define $T : \ell^n_p \longrightarrow 
L_\infty(\Cal A)$ by defining $T(e_i) = r_i$, $1 \leq i \leq n$ 
and extending linearly. Then $\|T\| = 1$.  Define $P : 
L_\infty(\Cal A) \longrightarrow \ell^n_p$ by $T(x) = 
\sum_{i=1}^n(\int xr_i)e_i$.  It is easily checked that $\|T : 
L_\infty(\Cal A) \longrightarrow \ell^n_1 \| \leq \sqrt{n}$; since 
$\| Id : \ell^n_1 \longrightarrow \ell^n_p \| = n^{\frac1p - 1}$, 
it follows that $\|P\| \leq n^{\frac1p - \frac12}$.  Obviously 
$PT = Id$, so up to a logarithmic factor, the order of 
$\underline\lambda(\ell^n_p)$ is correct. \endexample

\enddocument
\bye